\documentclass[10pt]{amsart}  
\usepackage{amsmath,amsthm,amssymb,amsfonts,graphicx, hyperref}

\title{Grothendieck inequalities for semidefinite programs with rank constraint}

\author{Jop Bri\"et}
\address{J.~Bri\"et, Centrum Wiskunde \& Informatica (CWI),
Science Park 123, 1098 SJ Amsterdam, The Netherlands}
\email{j.briet@cwi.nl}

\author{Fernando M\'ario de Oliveira Filho}
\address{F.M.~de Oliveira Filho, Institut f\"ur Mathematik,
Freie Universit\"at Berlin, Arnim\-allee 2, 14195 Berlin, Germany}
\email{fmario@mi.fu-berlin.de}

\author{Frank Vallentin} 
\address{F.~Vallentin, Delft Institute of Applied Mathematics, Technical University of Delft, P.O. Box 5031, 2600 GA Delft, The Netherlands}
\email{f.vallentin@tudelft.nl}

\thanks{The first author is supported by Vici grant 639.023.302 from the Netherlands Organization for Scientific Research (NWO), by the European Commission under the Integrated Project Qubit Applications (QAP) funded by the IST directorate as Contract Number 015848, by the Dutch BSIK/BRICKS project and by the European grant QCS. 
The second author was wupported by NWO Rubicon grant 680-50-1014.
The third author was supported by Vidi grant 639.032.917 from the Netherlands Organization for Scientific Research (NWO)}



\date{April 24, 20t12}

\newcommand{\R}{\mathbb{R}}

\newcommand{\Z}{\mathbb{Z}}
\newcommand{\C}{\mathbb{C}}
\newcommand{\E}{\mathbb{E}}

\newtheorem{defin}{Definition}[section]

\newtheorem{theorem}[defin]{Theorem}

\newtheorem{lemma}[defin]{Lemma}

\newtheorem{claim}{Claim}

\newenvironment{claimproof}[1][]{
  \begin{trivlist}
   \item[\hspace{\labelsep}{\sc\noindent Proof#1:\/}]}
   {{\hfill$\blacklozenge$}
  \end{trivlist}
}

\DeclareMathOperator{\sdp}{SDP}
\DeclareMathOperator{\sign}{sign}

\DeclareMathOperator{\rank}{rank}
\DeclareMathOperator{\arcsinh}{arcsinh}
\DeclareMathOperator{\Tr}{Tr}

\newcommand{\Hcal}{\mathcal{H}}
\newcommand{\bg}{\overline{g}}
\newcommand{\bE}{\overline{E}}

\newcommand{\beq}{\begin{equation}}
\newcommand{\eeq}{\end{equation}}
\newcommand{\beqn}{\begin{equation*}}
\newcommand{\eeqn}{\end{equation*}}
\newcommand{\beqr}{\begin{eqnarray}}
\newcommand{\eeqr}{\end{eqnarray}}
\newcommand{\beqrn}{\begin{eqnarray*}}
\newcommand{\eeqrn}{\end{eqnarray*}}

\newcommand{\sphere}[1]{S^{#1-1}} 
\newcommand{\HS}{\mathcal H}
\newcommand{\Exp}{\mathbb E}
\newcommand{\ind}{\mathbf{1}} 

\begin{document}

\begin{abstract}
Grothendieck inequalities are fundamental inequalities which are frequently used in many areas of mathematics and computer science. They can be interpreted as upper bounds for the integrality gap between two optimization problems: a difficult semidefinite program with rank-1 constraint and its easy semidefinite relaxation where the rank constrained is dropped. For instance, the integrality gap of the Goemans-Williamson approximation algorithm for MAX CUT can be seen as a Grothendieck inequality. In this paper we consider Grothendieck inequalities for ranks greater than $1$ and we give two applications: approximating ground states in the $n$-vector model in statistical mechanics and XOR games in quantum information theory.
\end{abstract}

\maketitle

\markboth{J.~Bri\"et, F.M.~de Oliveira Filho, F.~Vallentin}{Grothendieck inequalities for SDPs with rank constraint}

\section{Introduction}

Let $G = (V, E)$ be a graph with finite vertex set $V$ and edge set $E \subseteq \binom{V}{2}$. Let $A\colon V \times V \to \R$ be a symmetric matrix whose rows and columns are indexed by the vertex set of~$G$, and $r$ be a positive integer. The \emph{graphical Grothendieck problem with rank-$r$ constraint} is the following optimization problem:
\begin{equation*}
\sdp_r(G,A) = \max\biggl\{\,\sum_{\{u,v\} \in E} A(u,v)  f(u) \cdot f(v) \; : \; f\colon V \to S^{r-1}\,\biggr\},
\end{equation*}
where $S^{r-1} = \{\,x \in \R^r : x \cdot x = 1\,\}$ is the $(r-1)$-dimensional unit sphere. The \emph{rank-$r$ Grothendieck constant of the graph $G$} is the smallest constant $K(r,G)$ so that for all symmetric matrices $A\colon V \times V \to \R$ the following inequality holds:
\begin{equation}
\label{Grothendieck's inequality}
\sdp_{\infty}(G,A) \leq K(r,G) \sdp_r(G,A).
\end{equation}
Here $S^{\infty}$ denotes the unit sphere of the Hilbert space $l^2(\R)$ of square summable sequences, which contains $\R^n$ as the subspace of the first $n$ components. It is easy to see that $K(r,G) \geq 1$. 
In this paper, we prove new upper bounds for~$K(r,G)$.

\subsection{Some history}

Inequality \eqref{Grothendieck's inequality} is called a \emph{Grothendieck inequality} because it first appeared in the work \cite{Grothendieck} of Grothendieck on the metric theory of tensor products. More precisely, Grothendieck considered the case~$r = 1$ for $2$-chromatic (bipartite) graphs, although in quite a different language.  (A \emph{$k$-chromatic graph} is a graph whose chromatic number is $k$, i.e., one can color its vertices with $k$ colors so that adjacent vertices get different colors, but~$k-1$ colors do not suffice for this.) Grothendieck proved that in this case $K(1,G)$ is upper bounded by a constant that is independent of the size of $G$.

Later, Lindenstrauss and  Pe{\l}czy{\'n}ski \cite{LindenstraussPelczynski} reformulated Grothendieck's inequality for bipartite graphs in a way that is very close to the formulation we gave above. The graphical Grothendieck problem with rank-$1$ constraint was introduced by Alon, Makarychev, Makarychev, and Naor \cite{AlonMakarychevMakarychevNaor}. Haagerup \cite{Haagerup} considered the complex case of Grothendieck's inequality;
his upper bound is also valid for the real case~$r = 2$. The higher rank case for bipartite graphs was introduced by Bri\"et, Buhrman, and Toner~\cite{BrietBuhrmanToner}.

\subsection{Computational perspective}

There has been a recent surge of interest in Grothendieck inequalities by the computer science community. The problem $\sdp_r(G,A)$ is a semidefinite maximization problem with rank-$r$ constraint:
\begin{equation*}
\begin{split}
\sdp_r(G,A) = \max\biggl\{\,\sum_{\{u,v\}\in E} A(u,v)X(u,v) \;\; : \;\; & X \in \R^{V \times V}_{\succeq 0},\\[-1em] 
&  X(u,u) = 1 \text{ for all $u \in V$,}\\
& \rank X \leq r\,\biggr\},
\end{split}
\end{equation*}
where~$\R^{V \times V}_{\succeq 0}$ is the set of matrices~$X\colon V \times V \to \R$ that are positive semidefinite.

On the one hand,  $\sdp_r(G,A)$ is generally a difficult computational problem.  For instance, if $r=1$ and~$G$ is the complete bipartite graph $K_{n,n}$ on $2n$ nodes, and if~$A$ is the Laplacian matrix of a graph $G'$ on $n$ nodes, then computing $\sdp_1(K_{n,n},A)$ is equivalent to computing the weight of a maximum cut of~$G'$. The maximum cut problem (MAX CUT) is one of Karp's 21 $\mathrm{NP}$-complete problems.  On the other hand, if we relax the rank-$r$ constraint, then we deal with $\sdp_{\infty}(G,A)$, which is an easy computational problem: Obviously, one has $\sdp_{\infty}(G,A) = \sdp_{|V|}(G,A)$ and computing $\sdp_{|V|}(G,A)$ amounts to solving a semidefinite programming problem (see e.g.\ Vandenberghe, Boyd \cite{VandenbergheBoyd}). Therefore one may approximate it to any fixed precision in polynomial time by using the ellipsoid method or interior point algorithms. 

In many cases the optimal constant $K(r,G)$ is not known and so one is interested in finding upper bounds for~$K(r,G)$. Usually, proving an upper bound amounts to giving a randomized polynomial-time approximation algorithm for $\sdp_r(G,A)$.  In the case of the MAX CUT problem, Goemans and Williamson~\cite{GoemansWilliamson} pioneered an approach based on randomized rounding: One rounds an optimal solution of~$\sdp_\infty(G, A)$ to a feasible solution of~$\sdp_r(G, A)$. The expected value of the rounded solution is then related to the one of the original solution, and this gives an upper bound for~$K(r, G)$. Using this basic idea, Goemans and Williamson~\cite{GoemansWilliamson} showed that for all symmetric matrices $A\colon V \times V \to \R$ which have the properties $A(u,v) \leq 0$ for $u$ distinct from $v$ and $\sum_{u \in V} A(u,v) = 0$ for all $v \in V$, we have
\begin{equation*}
\sdp_\infty(K_{n,n},A) \leq (0.878\dots)^{-1} \sdp_1(K_{n,n},A).
\end{equation*}
\medskip

\subsection{Applications and references}

Grothendieck's inequality is a fundamental inequality in the theory of Banach spaces. Many books on the geometry of Banach spaces contain a substantial treatment of the result. We refer for instance to the books by Pisier~\cite{Pisier}, Jameson~\cite{Jameson}, and Garling~\cite{Garling}.

During the last years, especially after Alon and Naor \cite{AlonNaor} pointed out the connection between the inequality and approximation algorithms using semidefinite programs, Grothendieck's inequality has also become a unifying and fundamental tool outside of functional analysis. 

It has applications in optimization (Nesterov \cite{Nesterov}, Nemirovski, Roos, Terlaky \cite{NemirovskiRoosTerlaky}, Megretski \cite{Megretski}),  extremal combinatorics (Alon, Naor \cite{AlonNaor}), system theory (Ben-Tal, Nemirovski \cite{BenTalNemirovski}), machine learning (Charikar, Wirth \cite{CharikarWirth}, Khot, Naor \cite{KhotNaor, KhotNaor2}), communication complexity (Linial, Shraibman \cite{LinialSchraibman}), quantum information theory (Tsirel'son \cite{Tsirelson}, Regev, Toner \cite{RegevToner}), and computational complexity (Khot, O'Donnell \cite{KhotODonnell}, Arora, Berger, Kindler, Safra, Hazan \cite{AroraBergerKindlerHazanSafra}, Khot and Naor~\cite{KhotNaor3}, Raghavendra, Steurer \cite{RaghavendraSteurer}).

The references above mainly deal with the combinatorial rank $r = 1$ case, when $S^0 = \{-1,+1\}$. For applications in quantum information (Bri\"et, Buhrman, Toner \cite{BrietBuhrmanToner}) and in statistical mechanics (mentioned in Alon, Makarychev, Makarychev,  Naor \cite{AlonMakarychevMakarychevNaor}, Kindler, Naor, Schechtman \cite{KindlerNaorSchechtman}) the more geometrical case when $r > 1$ is of interest --- this case is the subject of this paper.

Before we present our results we consider the application to statistical mechanics: The \emph{$n$-vector model}, introduced by Stanley \cite{Stanley}, describes the interaction of particles in a spin glass with ferromagnetic and antiferromagnetic interactions. The case  $n = 1$ corresponds to the Ising model, the case $n = 2$ to the XY model, the case~$n = 3$ to the Heisenberg model, and the case $n = \infty$ to the Berlin-Kac spherical model. 

Let $G = (V, E)$ be the interaction graph where the vertices are particles and where edges indicate which particles interact.
The potential function $A\colon V \times V \to \R$ is $0$ if $u$ and $v$ are not adjacent, it is positive if there is ferromagnetic interaction between $u$ and $v$, and it is negative if there is antiferromagnetic interaction. The particles possess a vector-valued spin $f\colon V \to S^{n-1}$. In the absence of an external field, the total energy of the system is given by the \emph{Hamiltonian} 
\begin{equation*}
H(f) = -\sum_{\{u,v\} \in E} A(u,v)  f(u) \cdot f(v).
\end{equation*}
The ground state of this model is a configuration of spins $f\colon V \to S^{n-1}$ which minimizes the total energy. Finding the ground state is the same as solving~$\sdp_n(G,A)$. Typically, the interaction graph has small chromatic number, e.g.\ the most common case is when $G$ is a finite subgraph of the integer lattice $\Z^n$ where the vertices are the lattice points and where two vertices are connected if their Euclidean distance is one. These graphs are bipartite since they can be partitioned into even and odd vertices, corresponding to the parity of the sum of the coordinates.

We briefly describe the relation to quantum information theory. 
In an influential paper, Einstein, Podolsky, and Rosen~\cite{Einstein:1935} pointed out an anomaly of quantum mechanics that allows spatially separated parties to establish peculiar correlations by each performing measurements on a private quantum system: {\em entanglement}. 
Later, Bell~\cite{Bell:1964} proved that local measurements on a pair of spatially separated, entangled quantum systems, can give rise to joint probability distributions of measurement outcomes that violate certain inequalities (now called Bell inequalities), satisfied by any classical distribution.
Experimental results of Aspect, Grangier, and Roger~\cite{Aspect:1981} give strong evidence that nature indeed allows distant physical systems to be correlated in such non-classical ways.

{\em XOR games}, first formalized by Cleve, H\o yer, Toner, and Watrous~\cite{Cleve:2004}, constitute the simplest model in which entanglement can be studied quantitatively. In an XOR game, two players, Alice and Bob, receive questions $u$ and $v$ (resp.) that are picked by a referee according to some probability distribution $\pi(u,v)$ known to everybody in advance. Without sharing their questions, the players have to answer the referee with bits $a$ and $b$ (resp.), and win the game if  and only if the exclusive-OR of their answers $a\oplus b$ equals the value of a Boolean function $g(u,v)$; the function $g$ is also known in advance to all three parties.

In a quantum-mechanical setting, the players determine their answers by performing measurements on their shares of a pair of entangled quantum systems.
A {\em state} of a pair of $d$-dimensional quantum systems is a trace-$1$ positive semidefinite operator $\rho\in \C^{d^2\times d^2}_{\succeq 0}$.
The systems are {\em entangled} if $\rho$ is not a convex combination of tensor products of $d$-by-$d$ positive semidefinite matrices. 
For each question $u$, Alice has a two-outcome measurement defined by a pair of $d$-by-$d$ positive semidefinite matrices $\{A_u^0,A_u^1\}$ that satisfies $A_u^0 + A_u^1 = I$, where $I$ is the identity matrix.
Bob has a similar pair $\{B_v^0,B_v^1\}$ for each question~$v$. 
When the players perform their measurements, the probability that they obtain bits $a$ and $b$ is given by $\Tr(A_u^a\otimes B_v^b\rho)$.

The case $d = 1$ corresponds to a classical setting. In this case, the maximum winning probability equals $\big(1 + \sdp_1(G,A)\big)/2$, where
$G$ is the complete bipartite graph with Alice and Bob's questions on opposite sides of the partition, and $A(u,v) =  (-1)^{g(u,v)}\pi(u,v)/2$ for pairs $\{u,v\}\in E$ and $A(u,v)=0$ everywhere else.

Tsirel'son~\cite{Tsirelson} related the maximum winning probability $\omega^*_d(\pi,g)$ of the game $(\pi,g)$, when the players are restricted to measurements on $d$-dimensional quantum systems, to the quantity $\sdp_r(G,A)$.
In particular, he proved that
\begin{equation*}
\frac{1 + \sdp_{\lfloor\log d\rfloor}(G,A)}{2}  \leq \omega_d^*(\pi,g) \leq \frac{1+ \sdp_{2d}(G,A)}{2}.
\end{equation*}
The quantity~$\sdp_r(G,A)$ thus gives bounds on the maximum winning probability of XOR games when players are limited in the amount of entanglement they are allowed to use. 
The rank-$r$ Grothendieck constant $K(r,G)$ of the bipartite graph~$G$ described above gives a quantitative bound on the advantage that unbounded entanglement gives over finite entanglement in XOR games.
\medskip
\medskip
\medskip

\subsection{Our results and methods}
\label{sec:our methods}

The purpose of this paper is to prove explicit upper bounds for $K(r,G)$. We are especially interested in the case of small~$r$ and graphs with small chromatic number, although our methods are not restricted to this. The proof of the following theorem gives a randomized polynomial-time approximation algorithm for approximating ground states in the Heisenberg model in the lattice $\Z^3$ with approximation ratio $0.78\ldots = (1.28\ldots)^{-1}$. This result can be regarded as one of the main contributions of this paper.
\newpage
\begin{theorem}
\label{thm:main}
For $r = 1, \ldots, 10$ and in the case of a bipartite or a tripartite graph~$G$ the rank-$r$ Grothendieck constant is at most:
\medskip

\begin{center}
\begin{tabular}{ccc}
\noalign{\hrule\vskip2pt}
\quad$r$\quad\hbox{} & \quad{\sl bipartite $G$}\quad \hbox{} & \quad {\sl tripartite $G$} \quad\hbox{}\\[2pt]
\noalign{\hrule\vskip2pt}
$1$ & $1.782213\dots$ & $3.264251\dots$\\
$2$  &                              $1.404909\dots$  & $2.621596\dots$\\
$3$  &                              $1.280812\dots$  & $2.412700\dots$\\
$4$  &                              $1.216786\dots$  & $2.309224\dots$\\
$5$  &                              $1.177179\dots$  & $2.247399\dots$\\
$6$  &                              $1.150060\dots$  & $2.206258\dots$\\
$7$  &                              $1.130249\dots$  & $2.176891\dots$\\
$8$  &                              $1.115110\dots$  & $2.154868\dots$\\
$9$  &                              $1.103150\dots$  & $2.137736\dots$\\
$10$ &                             $1.093456\dots$  & $2.124024\dots$\\[2pt]
\hline
\end{tabular}
\end{center}
\end{theorem}
\smallskip

The bound for the original Grothendieck constant $K(1,G)$ for bipartite~$G$ is due to Krivine~\cite{Krivine}. 
For more than thirty years this was the best known upper bound, and it was conjectured by many to be optimal.
However, shortly after our work appeared in preprint form, Braverman, Makarychev, Makarychev and Naor~\cite{Braverman:2011} showed that  Krivine's bound can be slightly improved.
The best known lower bound is $1.676956\dots$ due to Davie~\cite{Davie} and Reeds~\cite{Reeds} (see also Khot and O'Donnell~\cite{KhotODonnell}). The bound for~$K(2,G)$ is due to Haagerup \cite{Haagerup}. 

When the graph~$G$ has large chromatic number, then the result of Alon, Ma\-ka\-ry\-chev, Makarychev, and Naor~\cite{AlonMakarychevMakarychevNaor} gives the best known bounds for $K(1,G)$: They prove a logarithmic dependence on the chromatic number of the graph (actually on the theta number of the complement of $G$, cf. Section~\ref{constructing section}) whereas our methods only give a linear dependence.  Although our main focus is on small chromatic numbers, for completeness we extend the results of~\cite{AlonMakarychevMakarychevNaor} for large chromatic numbers to $r\geq 2$ in Section~\ref{sec:highchrom}. In a previous paper~\cite{Briet:MTNS2010} we proved that $K(r,K_{n,n}) = 1 +\Theta(1/r)$.

For the proof of Theorem~\ref{thm:main} we use the framework of Krivine and Haagerup which we explain below. Our main technical contributions are a matrix version of Grothendieck's identity (Lemma~\ref{grothendieck identity}) and a method to construct new unit vectors which can also deal with nonbipartite graphs (Lemma~\ref{lem:gen-embedding}).

The strategy of Haagerup and Krivine is based on the following embedding lemma:

\begin{lemma}
\label{lem:embed-vague}
Let~$G = (V, E)$ be a graph and choose~$Z = (Z_{ij}) \in \R^{r \times |V|}$ at random so that each entry is distributed independently according to the normal distribution with mean~$0$ and variance~$1$, that is,~$Z_{ij} \sim N(0, 1)$.

Given~$f\colon V \to S^{|V|-1}$, there is a function~$g\colon V \to S^{|V|-1}$ such that whenever~$u$ and~$v$ are adjacent in~$G$, then
\[
\E\biggl[\frac{Z g(u)}{\|Z g(u)\|} \cdot \frac{Z g(v)}{\|Z g(v)\|}\biggr] = \beta(r, G) f(u) \cdot f(v)
\]
for some constant~$\beta(r, G)$ depending only on~$r$ and~$G$.
\end{lemma}

In the statement above we are vague regarding the constant~$\beta(r, G)$. We will give the precise statement of the lemma in Section~\ref{constructing section} (cf.~Lemma~\ref{lem:gen-embedding} there), right now this precise statement is not relevant to our discussion.

Now, the strategy of Haagerup and Krivine amounts to analyzing the following four-step procedure that yields  a randomized polynomial-time approximation algorithm for $\sdp_r(G,A)$:

\medskip

\noindent
{\bf Algorithm~A.}\enspace Takes as input a finite graph~$G = (V, E)$ with at least one edge and a symmetric matrix~$A\colon V \times V \to \R$, and returns a feasible solution~$h\colon V \to S^{r-1}$ of~$\sdp_r(G, A)$.

\begin{enumerate}
\item Solve $\sdp_{\infty}(G, A)$, obtaining an optimal solution $f\colon V \to S^{|V|-1}$.
\item Use $f$ to construct $g\colon V \to S^{|V|-1}$ according to Lemma~\ref{lem:embed-vague}.
\item Choose $Z = (Z_{ij}) \in \R^{r \times |V|}$ at random so that every matrix entry $Z_{ij}$ is distributed independently according to the standard normal distribution with mean $0$ and variance $1$, that is,~$Z_{ij} \sim N(0,1)$.
\item Define $h\colon V \to S^{r-1}$ by setting $h(u) = Zg(u)/\|Zg(u)\|$.
\end{enumerate}
\medskip

To analyze this procedure, we compute the expected value of the feasible solution~$h$. Using Lemma~\ref{lem:embed-vague} we obtain
\begin{equation}
\label{ineqchain}
\begin{split}
\sdp_r(G,A)  & \geq \E\biggl[\sum_{\{u,v\} \in E} A(u,v) h(u) \cdot h(v)\biggr]\\
&  =  \sum_{\{u,v\} \in E} A(u,v) \E[h(u) \cdot h(v)]\\
&  =  \beta(r,G) \sum_{\{u,v\} \in E} A(u,v) f(u) \cdot f(v)\\
&  =  \beta(r,G) \sdp_{\infty}(G,A),
\end{split}
\end{equation}
and so we have $K(r,G) \leq \beta(r,G)^{-1}$. 

If we were to skip step~{(2)} and apply step~{(4)} to $f$ directly, then the expectation~$\E[h(u)\cdot h(v)]$ would be a non-linear function of $f(u)\cdot f(v)$, which would make it difficult to assess the quality of the feasible solution $h$.  The purpose of step~{(2)} is to linearize this expectation, which allows us to estimate the quality of~$h$ in terms of a linear function of $\sdp_r(G,A)$.

The constant~$\beta(r, G)$ in Lemma~\ref{lem:embed-vague} is defined in terms of the Taylor expansion of the inverse of the function~$E_r\colon [-1, 1] \to [-1, 1]$ given by
\[
E_r(x \cdot y) = \E\biggl[\frac{Zx}{\|Zx\|} \cdot \frac{Zy}{\|Zy\|}\biggr],
\]
where~$x$, $y \in S^\infty$ and~$Z = (Z_{ij}) \in \R^{r \times \infty}$ is chosen so that its entries are independently distributed according to the normal distribution with mean~$0$ and variance~$1$. The function~$E_r$ is well-defined since the expectation above is invariant under orthogonal transformations.

The Taylor expansion of~$E_r$ is computed in Section~\ref{grothendieck identity section}. The Taylor expansion of~$E_r^{-1}$ is treated in Section~\ref{convergence section}, where we basically follow Haagerup~\cite{Haagerup}. A precise version of Lemma~\ref{lem:embed-vague} is stated and proved in Section~\ref{constructing section}, following Krivine~\cite{Krivine}. 

Finally, in Section~\ref{refined section} we show that one can refine this analysis and can (strictly) improve the upper bound if one takes the dimension of the matrix $A\colon V \times V \to \R$ into account. In particular, we compare the problems $\sdp_q$ and $\sdp_r$ for $q \geq r$. Earlier, Avidor and Zwick~\cite{AvidorZwick} considered the problem of bounding the ratio $\sdp_q(G,A)/\sdp_1(G,A)$ for $q = 2,3$ and $A$ the Laplacian matrix of a graph.

\section{A matrix version of Grothendieck's identity}
\label{grothendieck identity section}

In the analysis of many approximation algorithms that use semidefinite programming the following identity plays a central role:
Let $u$, $v$ be unit (column) vectors in $\R^n$ and let $Z \in \R^{1 \times n}$ be a random (row) vector whose entries are distributed independently according to the standard normal distribution with mean $0$ and variance~$1$. Then,
\begin{equation*}
\E[\sign(Zu)\sign(Zv)] = \E\biggl[\frac{Zu}{\|Zu\|} \cdot \frac{Zv}{\|Zv\|}\biggr] = \frac{2}{\pi}\arcsin(u \cdot v).
\end{equation*}

For instance, the celebrated algorithm of Goemans and Williamson \cite{GoemansWilliamson} for approximating the MAX CUT problem is based on this. The identity is called \emph{Grothendieck's identity} since it appeared for the first time in Grothendieck's work on the metric theory of tensor products \cite[Proposition 4, p. 63]{Grothendieck} (see also Diestel, Fourie, and Swart~\cite{DiestelFourieSwart}).

In this section we extend Grothendieck's identity from vectors to matrices by replacing the arcsine function by a hypergeometric function.

\begin{lemma}
\label{grothendieck identity}
Let $u$, $v$ be unit vectors in $\R^n$ and let $Z \in \R^{r \times n}$  be a random matrix whose entries are distributed independently according to the standard normal distribution with mean $0$ and variance $1$. Then,
\begin{equation*}
\E\biggl[\frac{Zu}{\|Zu\|} \cdot \frac{Zv}{\|Zv\|}\biggr]= \frac{2}{r}\left(\frac{\Gamma((r+1)/2)}{\Gamma(r/2)}\right)^2 
(u \cdot v)\, {}_{2}F_{1} \left(\!\!\begin{array}{cc} 1/2, 1/2\\ r/2+1\end{array}\!; (u \cdot v)^2 \right).
\end{equation*}
Here,
\begin{equation*}
\begin{split}
&(u \cdot v)\, {}_{2}F_{1} \left(\!\!\begin{array}{cc} 1/2, 1/2\\ r/2+1\end{array}\!; (u \cdot v)^2 \right)\\
&\qquad{}= \sum_{k = 0}^{\infty} \frac{(1\cdot 3 \cdots (2k-1))^2}{(2 \cdot 4 \cdots 2k)((r+2)\cdot(r+4)\cdots (r+2k))}  (u \cdot v)^{2k+1}.
\end{split}
\end{equation*}
\end{lemma}

Before proving the lemma we review special cases known in the literature. If $r = 1$, then we get the original Grothendieck's identity:
\begin{eqnarray*}
\E[\sign(Zu)\sign(Zv)] 
& = & \frac{2}{\pi} \arcsin (u \cdot v)\\
&  = & \frac{2}{\pi} \left( u \cdot v + \left(\frac{1}{2}\right)\frac{(u \cdot v)^3}{3} + \left(\frac{1\cdot 3}{2 \cdot 4}\right) \frac{(u \cdot v)^5}{5} + \cdots \right).
\end{eqnarray*}
The case $r=2$ is due to Haagerup \cite{Haagerup}: 
\begin{eqnarray*}
\E\left[\frac{Zu}{\|Zu\|} \cdot \frac{Zv}{\|Zv\|}\right]
& = & \frac{1}{u \cdot v}\left(E(u \cdot v)-(1-(u \cdot v)^2)K(u \cdot v)\right)\\
& = & \frac{\pi}{4}
\left( 
u \cdot v + \left(\frac{1}{2}\right)^2 \frac{(u \cdot v)^3}{2} 
+ \left(\frac{1 \cdot 3}{2 \cdot 4}\right)^2 \frac{(u \cdot v) ^5}{3}  + \cdots\right),
\end{eqnarray*}
where $K$ and $E$ are the complete elliptic integrals of the first and second kind. Note that on page 201 of Haagerup~\cite{Haagerup} $\pi/2$ has to be $\pi/4$. Bri\"et, Oliveira, and Vallentin~\cite{BrietOliveiraVallentin}  computed the first coefficient $2/r(\Gamma((r+1)/2)/\Gamma(r/2))^2$ of the Taylor series of the expectation for every $r$.

The following elegant proof of Grothendieck's identity has become a classic: We have $\sign(Zu)\sign(Zv) = 1$ if and only if the vectors $u$ and $v$ lie on the same side of the hyperplane orthogonal to the vector $Z \in \R^{1 \times n}$. Now we project this $n$-dimensional situation to the plane spanned by $u$ and $v$. Then the projected random hyperplane becomes a random line. This random line is distributed according to the uniform probability measure on the unit circle because $Z$ is normally distributed. Now one obtains the final result by measuring intervals on the unit circle: The probability that $u$ and $v$ lie on the same side of the line is $1 - \arccos(u \cdot v)/\pi$.

However, we do not have such a picture proof for our matrix version. Our proof is based on the rotational invariance of the normal distribution and integration with respect to spherical coordinates together with some identities for hypergeometric functions. A similar calculation was done by K\"onig and Tomczak-Jaegermann~\cite{Koenig}. It would be interesting to find a more geometrical proof of the lemma.

For computing the first coefficient of the Taylor series in \cite{BrietOliveiraVallentin} we took a slightly different route: We integrated using the Wishart distribution of $2 \times 2$-matrices.

\begin{proof}[Proof of Lemma \ref{grothendieck identity}]
Let $Z_i \in \R^n$ be the $i$-th row of the matrix $Z$, with $i = 1, \ldots r$. We define vectors
\begin{equation*}
x =
\begin{pmatrix}
Z_1 \cdot u\\
Z_2 \cdot u\\
\vdots\\
Z_r \cdot u
\end{pmatrix}\qquad\text{and}\qquad
y = 
\begin{pmatrix}
Z_1 \cdot v\\
Z_2 \cdot v\\
\vdots\\
Z_r \cdot v
\end{pmatrix}
\end{equation*}
so that we have $x\cdot y = (Zu)\cdot (Zv)$.
Since the probability distribution of the vectors $Z_i$ is invariant under orthogonal transformations we may assume that $u = (1, 0, \ldots, 0)$ and $v = (t, \sqrt{1-t^2}, 0, \ldots, 0)$ and so the pair $(x,y) \in \R^r \times \R^r$ is distributed according to the probability density function~(see e.g.\ Muirhead~\cite[p.~10, eq.~(7)]{Muirhead})
\begin{equation*}
(2\pi \sqrt{1-t^2})^{-r}\exp\left(-\frac{x \cdot x - 2t x \cdot y + y \cdot y}{2(1-t^2)}\right).
\end{equation*}
Hence,
\begin{equation*}
\begin{split}
&\E\left[\frac{x}{\|x\|} \cdot \frac{y}{\|y\|}\right]\\
&\qquad{}= (2\pi \sqrt{1-t^2})^{-r} \int_{\R^r} \int_{\R^r} \frac{x}{\|x\|} \cdot \frac{y}{\|y\|} \exp\left(-\frac{x \cdot x - 2t x \cdot y + y \cdot y}{2(1-t^2)}\right)\, dx dy.
\end{split}
\end{equation*}
By using spherical coordinates $x = \alpha \xi$, $y = \beta \eta$, where $\alpha,\beta \in [0,\infty)$ and $\xi, \eta \in S^{r-1}$, we rewrite the above integral as
\begin{equation*}
\vcenter{\halign{\hbox to\hsize{#\hfil}\cr
\quad $\displaystyle\int_0^{\infty} \int_0^{\infty} (\alpha\beta)^{r-1} \exp\left(-\frac{\alpha^2+\beta^2}{2(1-t^2)}\right)\int_{S^{r-1}} \int_{S^{r-1}} \xi \cdot \eta \exp\left(\frac{\alpha\beta t\xi \cdot \eta}{1-t^2}\right)$\hfill\cr\noalign{\vskip1pt}
\hfill $d\omega(\xi) d\omega(\eta) d\alpha d\beta.$\quad\cr
}}
\end{equation*}

If $r = 1$, we get for the inner double integral
\begin{equation*}
\begin{split}
& \int_{S^{0}} \int_{S^{0}} \xi \cdot \eta \exp\left(\frac{\alpha\beta t\xi \cdot \eta}{1-t^2}\right)\, d\omega(\xi) d\omega(\eta)\\
&\qquad{}= 4 \sinh\left(\frac{\alpha\beta t}{1-t^2}\right)\\
&\qquad{}= 4 \frac{\alpha\beta t}{1-t^2} {}_{0}F_{1} \left(\!\!\begin{array}{cc} \overline{\phantom{xxx}}\\ 3/2\end{array}\!; \left(\frac{\alpha\beta t}{2(1-t^2)}\right)^2 \right).
\end{split}
\end{equation*}

Now we consider the case when $r \geq 2$. Since the inner double integral over the spheres only depends on the inner product $p= \xi \cdot \eta$ it can be rewritten as
\begin{equation*}
\omega(S^{r-2}) \omega(S^{r-1})\int_{-1}^1 p \exp\left(\frac{\alpha\beta tp}{1-t^2}\right) (1-p^2)^{(r-3)/2}\, dp,
\end{equation*}
where
\begin{equation*}
\omega(S^{r-2}) \omega(S^{r-1}) = \frac{4 \pi^{r-1/2}}{\Gamma(r/2)\Gamma((r-1)/2)}.
\end{equation*}
Integration by parts yields
\begin{equation*}
\begin{split}
& \int_{-1}^1 p(1-p^2)^{(r-3)/2} \exp\left(\frac{\alpha\beta tp}{1-t^2}\right)\, dp\\
&\qquad{}= \frac{\alpha\beta t}{(r-1)(1-t^2)} \int_{-1}^1 (1-p^2)^{(r-1)/2} \exp\left(\frac{\alpha\beta tp}{1-t^2}\right)  dp.
\end{split}
\end{equation*}
The last integral can be rewritten using the modified Bessel function of the first kind (cf.~Andrews, Askey, Roy \cite[p.~235, Exercise 9]{AndrewsAskeyRoy})
\begin{equation*}
\begin{split}
& \int_{-1}^1 (1-p^2)^{(r-1)/2} \exp\left(\frac{\alpha\beta tp}{1-t^2}\right)\, dp\\
&\qquad{}= \Gamma((r+1)/2) \sqrt{\pi} \left(\frac{2(1-t^2)}{\alpha\beta t}\right)^{r/2} I_{r/2}\left(\frac{\alpha\beta t}{1-t^2}\right).
\end{split}
\end{equation*}
One can write $I_{r/2}$ as a hypergeometric function (cf.~Andrews, Askey, and Roy~\cite[(4.12.2)]{AndrewsAskeyRoy})
\begin{equation*}
I_{r/2}(x) = (x/2)^{r/2} \sum_{k=0}^{\infty} \frac{(x/2)^{2k}}{k!\Gamma(r/2+k+1)} = 
\frac{(x/2)^{r/2}}{\Gamma((r+2)/2)} {}_{0}F_{1} \left(\!\!\begin{array}{cc} \overline{\phantom{xxx}}\\ (r+2)/2\end{array}\!; \left(\frac{x}{2}\right)^2 \right).
\end{equation*}
Putting things together, we get 
\begin{equation*}
\begin{split}
& \omega(S^{r-2}) \omega(S^{r-1})\int_{-1}^1 p \exp\left(\frac{\alpha\beta tp}{1-t^2}\right) (1-p^2)^{(r-3)/2}\, dp\\
&\qquad{}= \frac{4\pi^r}{\Gamma(r/2)^2 r}
\frac{\alpha\beta t}{1-t^2} 
{}_{0}F_{1} \left(\!\!\begin{array}{cc} \overline{\phantom{xxx}}\\ (r+2)/2\end{array}\!; \left(\frac{\alpha\beta t}{2(1-t^2)}\right)^2 \right).
\end{split}
\end{equation*}
Notice that the last formula also holds for $r = 1$. So we can continue without case distinction.

Now we evaluate the outer double integral
\begin{equation*}
\int_0^{\infty}\int_0^{\infty} (\alpha\beta)^r \exp\left(-\frac{\alpha^2+\beta^2}{2(1-t^2)}\right) {}_{0}F_{1} \left(\!\!\begin{array}{cc} \overline{\phantom{xxx}} \\ (r+2)/2\end{array}\!; \left(\frac{\alpha\beta t}{2(1-t^2)}\right)^2 \right)\, d\alpha d\beta.
\end{equation*}
Here the inner integral equals
\begin{equation*}
\int_0^{\infty} \alpha^r \exp\left(-\frac{\alpha^2}{2(1-t^2)}\right) {}_{0}F_{1} \left(\!\!\begin{array}{cc} \overline{\phantom{xxx}} \\ (r+2)/2\end{array}\!; \left(\frac{\alpha\beta t}{2(1-t^2)}\right)^2 \right)\, d\alpha,
\end{equation*}
and doing the substitution $\gamma = \alpha^2 / (2(1-t^2))$ gives
\begin{equation*}
2^{(r-1)/2}(1-t^2)^{(r+1)/2}\int_0^{\infty} \gamma^{(r-1)/2} \exp(-\gamma)\, {}_{0}F_{1} \left(\!\!\begin{array}{cc} \overline{\phantom{xxx}} \\ (r+2)/2\end{array}\!; \frac{\gamma(\beta t)^2}{2(1-t^2)} \right)\, d\gamma,
\end{equation*}
which is by the Bateman Manuscript Project \cite [p.~337 (11)]{ErdelyiMagnusOberhettinerTricomi} equal to
\begin{equation*}
2^{(r-1)/2}(1-t^2)^{(r+1)/2} \Gamma((r+1)/2) {}_{1}F_{1} \left(\!\!\begin{array}{cc} (r+1)/2 \\ (r+2)/2\end{array}\!; \frac{(\beta t)^2}{2(1-t^2)} \right).
\end{equation*}
Now we treat the remaining outer integral in a similar way, using \cite [p.~219 (17)]{ErdelyiMagnusOberhettinerTricomi}, and get that
\begin{equation*}
\begin{split}
&\int_0^{\infty} \beta^r \exp\left(-\frac{\beta^2}{2(1-t^2)}\right) {}_{1}F_{1} \left(\!\!\begin{array}{cc} (r+1)/2 \\ (r+2)/2\end{array}\!; \frac{(\beta t)^2}{2(1-t^2)} \right)\, d\beta\\
&\qquad{}=2^{(r-1)/2}(1-t^2)^{(r+1)/2} \Gamma((r+1)/2) {}_{2}F_{1} \left(\!\!\begin{array}{cc}(r+1)/2, (r+1)/2 \\ (r+2)/2\end{array}\!; t^2 \right).
\end{split}
\end{equation*}

By applying Euler's transformation (cf.~Andrews, Askey, and Roy~\cite[(2.2.7)]{AndrewsAskeyRoy})  
\begin{equation*}
{}_{2}F_{1} \left(\!\!\begin{array}{cc}(r+1)/2, (r+1)/2 \\ (r+2)/2\end{array}\!; t^2 \right) = (1-t^2)^{-r/2} {}_{2}F_{1} \left(\!\!\begin{array}{cc}1/2, 1/2 \\ (r+2)/2\end{array}\!; t^2 \right)
\end{equation*}
and after collecting the remaining factors we arrive at the result.
\end{proof}

\section{Convergence radius}
\label{convergence section}

To construct the new vectors in the third step of the algorithm that are used to linearize the expectation we will make use of the Taylor series expansion of the inverse of~$E_r$. Locally around zero we can expand the function $E_r^{-1}$ as
\begin{equation*}
E_r^{-1}(t) = \sum_{k = 0}^{\infty} b_{2k+1} t^{2k+1},
\end{equation*}
but in the proof of Lemma~\ref{lem:gen-embedding} it will be essential that this expansion be valid for all $t \in [-1,1]$.

In the case $r = 1$ we have $E_1^{-1}(t) = \sin(\pi/2 t)$ and here the convergence radius is even infinity. The case $r = 2$ was treated by Haagerup and it requires quite some technical work which we sketch very briefly now. He shows that $|b_{k}| \leq C/k^2$ for some constant $C$, independent of $k$, using tools from complex analysis. Using Cauchy's integral formula and after doing some simplifications \cite[p.~208]{Haagerup} one can express $b_k$ as 
\begin{equation*}
b_k = \frac{2}{\pi k} \int_1^{\alpha} \Im(E_2(z)^{-k})\, dz + \frac{2}{\pi k}\Im\left(\int_{C'_{\alpha}} E_2(z)^{-k}\, dz\right),
\end{equation*}
where $C'_{\alpha}$ is the quarter circle $\{\,\alpha e^{i\theta} : \theta \in [0,\pi/2]\,\}$.

For an appropriate choice of $\alpha$ the first integral is in absolute value bounded above by $C/k$ and the second integral is in absolute value exponentially small in~$k$. We refer to the original paper for the details. One key point in the arguments is the following integral representation of $E_2$ giving an analytic continuation of $E_2$ on the complex plane slit along the half line $(1, \infty)$:
\begin{equation*}
E_2(z) = \int_0^{\pi/2} \sin\theta \arcsin(z \sin\theta)\, d\theta.
\end{equation*}
Here, the term $\arcsin(z \sin \theta)$ gives the main contribution in the estimates. 

Now we derive a similar representation of $E_r$ and using it in Haagerup's analysis with obvious changes shows that also for $r > 2$ we have $b_{k} \leq C/k^2$ for some constant $C$, independent of $k$.

\begin{lemma}
For $r \geq 2$ we have
\begin{equation*}
E_r(z) = \frac{2(r-1)\Gamma((r+1)/2)}{\Gamma(1/2)\Gamma(r/2)} \int_0^{\pi/2} \cos^{r-2} \theta \sin\theta \arcsin(z \sin \theta)\, d\theta.
\end{equation*}
\end{lemma}

\begin{proof}
Using Euler's integral representation of the hypergeometric function (cf.~Andrews, Askey, and Roy~\cite[Theorem 2.2.1]{AndrewsAskeyRoy}) we can rewrite~$E_r$ as
\begin{equation*}
E_r(z) = \frac{\Gamma((r+1)/2)}{\Gamma(1/2)\Gamma(r/2)}\int_0^1 \frac{(1-t)^{(r-1)/2}z}{\sqrt{t(1-z^2t)}}\, dt,
\end{equation*}
which is valid in the complex plane slit along the half line $(1,\infty)$. Using the substitution $t = \sin^2 \theta$ we get
\begin{equation*}
E_r(z) = 2 \frac{\Gamma((r+1)/2)}{\Gamma(1/2)\Gamma(r/2)} \int_0^{\pi/2} \frac{\cos^r \theta z}{\sqrt{1 - z^2 \sin^2 \theta}}\, d\theta.
\end{equation*}
Now integration by parts and the identity  
\begin{equation*}
\frac{d}{d\theta} \arcsin(z \sin \theta) = \frac{z \cos\theta}{\sqrt{1-z^2 \sin^2 \theta}}
\end{equation*}
gives the result.
\end{proof}

\section{Constructing new vectors}
\label{constructing section}

In this section we use the Taylor expansion of the inverse of the function~$E_r$ to give a precise statement and proof of Lemma~\ref{lem:embed-vague}; this is done in Lemma~\ref{lem:gen-embedding}. For this we follow Krivine~\cite{Krivine}, who proved the statement of the lemma in the case of bipartite graphs. We comment on how his ideas are related to our construction, which can also deal with nonbipartite graphs, after we prove the lemma.

For the nonbipartite case we need to use the theta number, which is a graph parameter introduced by  Lov\'asz~\cite{Lovasz}. Let~$G = (V, E)$ be a graph. The \emph{theta number} of the complement of~$G$, denoted by~$\vartheta(\overline{G})$, is the optimal value of the following semidefinite program:
\begin{equation}
\label{opt:theta-gbar}
\begin{split}
\vartheta(\overline{G}) = \min\Big\{\,\lambda : \;\; & Z \in \R^{V \times V}_{\succeq 0},\\
&Z(u, u) = \lambda - 1 \;\;\text{for~$u \in V$},\\
&Z(u, v) = -1 \;\; \text{for~$\{u,v\} \in E$}\,\Big\}.
\end{split}
\end{equation}
It is known that the theta number of the complement of~$G$ provides a lower bound for the chromatic number of~$G$. This can be easily seen as follows. Any proper $k$-coloring of $G$ defines a mapping of~$V$ to the vertices of a $(k-1)$-dimensional regular simplex whose vertices lie in a sphere of radius $\sqrt{k-1}$: Vertices in the graph having the same color are sent to the same vertex in the regular simplex and vertices of different colors are sent to different vertices in the regular simplex. The Gram matrix of these vectors gives a feasible solution of~\eqref{opt:theta-gbar}. 

\begin{lemma}
\label{lem:gen-embedding}
Let~$G = (V, E)$ be a graph with at least one edge. Given~$f\colon V
\to S^{|V|-1}$, there is a function~$g\colon V \to S^{|V|-1}$ such
that whenever $u$ and $v$ are adjacent, then
\begin{equation*}
E_r\big(g(u) \cdot g(v)\big) = \beta(r,G) f(u) \cdot f(v).
\end{equation*}
The constant $\beta(r,G)$ is defined as the solution of the equation
\begin{equation*}
\sum_{k=0}^{\infty} |b_{2k+1}| \beta(r,G)^{2k+1} = \frac{1}{\vartheta(\overline{G}) - 1},
\end{equation*}
where
\begin{equation*}
E_r^{-1}(t) = \sum_{k=0}^{\infty} b_{2k+1} t^{2k+1}.
\end{equation*}
\end{lemma}

With this lemma, we can give a proof of Theorem~\ref{thm:main}.

\begin{proof}[Proof of Theorem~\ref{thm:main}]
We combine Lemma~\ref{lem:gen-embedding} with the analysis of Algorithm~A from Section~\ref{sec:our methods}. To compute the table in the theorem, we use the formula
\begin{equation}\label{eq:rankgroth-mfinverse}
b_k = \frac{1}{k!a_1^k} \left[\frac{d^{k-1}}{dt^{k-1}}\left(1 + \frac{a_2}{a_1}t + \cdots + \frac{a_k}{a_1}t^{k-1}\right)^{-k}\right]_{t = 0},
\end{equation}
where $a_i$ are the Taylor coefficients of $E_r$ (cf.~Morse and Feshbach~\cite[(4.5.13)]{MorseFeshbach}).
\end{proof}

Now we give a proof of the lemma.

\begin{proof}[Proof of Lemma~\ref{lem:gen-embedding}]
We construct the vectors $g(u) \in S^{|V|-1}$ by constructing vectors $R(u)$ in an infinite-dimensional Hilbert space whose inner product matrix coincides with the one of the $g(u)$. We do this in three steps.

In the first step, set $H = \R^{|V|}$ and consider the Hilbert space
\begin{equation*}
\mathcal{H} = \bigoplus_{k=0}^\infty H^{\otimes (2k+1)}.
\end{equation*}
For a unit vector~$x \in H$, consider the vectors~$S(x)$, $T(x) \in \mathcal{H}$ given componentwise by
\begin{equation*}
S(x)_k = \sqrt{|b_{2k+1}| \beta(r,G)^{2k+1}} x^{\otimes (2k+1)}
\end{equation*}
and
\begin{equation*}
T(x)_k = \sign(b_{2k+1}) \sqrt{|b_{2k+1}| \beta(r,G)^{2k+1}} x^{\otimes (2k+1)}.
\end{equation*}
Then for vectors $x, y \in S^{|V|-1}$ we have
\begin{equation*}
S(x) \cdot T(y) = E_r^{-1}(\beta(r,G) x \cdot y)
\end{equation*}
and moreover
\begin{equation*}
S(x) \cdot S(x) = T(x) \cdot T(x) = \sum_{k=0}^{\infty} |b_{2k+1}| \beta(r,G)^{2k+1} = \frac{1}{\vartheta(\overline{G}) - 1}.
\end{equation*}

In the second step, let ~$\lambda = \vartheta(\overline{G})$, and $Z$ be an optimal solution of~\eqref{opt:theta-gbar}. We have~$\lambda \geq 2$ since $G$ has at least one edge. Set
\begin{equation*}
A = \frac{(\lambda - 1) (J+Z)}{2\lambda}\qquad
\text{and}\qquad
B = \frac{(\lambda - 1) J - Z}{2\lambda},
\end{equation*}
and consider the matrix
\begin{equation*}
U = \begin{pmatrix}
A&B\\
B&A
\end{pmatrix}.
\end{equation*}
By applying a Hadamard transformation
\begin{equation*}
\frac{1}{\sqrt{2}} 
\begin{pmatrix}
I & I\\
I & -I
\end{pmatrix}
U 
\frac{1}{\sqrt{2}} 
\begin{pmatrix}
I & I\\
I & -I
\end{pmatrix}
=
\begin{pmatrix}
A+B & 0\\
0 & A - B
\end{pmatrix}
\end{equation*}
one sees that $U$ is positive semidefinite, since both~$A + B$ and~$A - B$ are positive semidefinite. Define $s\colon V \to \R^{2|V|}$ and $t\colon V \to \R^{2|V|}$ so that
\begin{equation*}
s(u) \cdot s(v) = t(u) \cdot t(v) = A(u,v)\qquad\text{and}\qquad s(u) \cdot t(v) = B(u,v).
\end{equation*}
The matrix $U$ is the Gram matrix of the vectors $\big(s(u)\big)_{u\in V}$ and $\big(t(v)\big)_{v\in V}$. 
It follows that these maps have the following properties:
\begin{enumerate}
\item $s(u) \cdot t(u) = 0$ for all~$u \in V$,
\item $s(u) \cdot s(u) = t(u) \cdot t(u) = (\vartheta(\overline{G}) - 1)/2$ for all~$u \in V$,
\item $s(u) \cdot s(v) = t(u) \cdot t(v) = 0$ whenever~$\{u,v\} \in E$,
\item $s(u) \cdot t(v) = s(v) \cdot t(u) = 1/2$ whenever~$\{u,v\} \in E$.
\end{enumerate}

In the third step we combine the previous two. We define the vectors
\begin{equation*}
R(u) = s(u) \otimes S(f(u)) + t(u) \otimes T(f(u)).
\end{equation*}
For adjacent vertices $u, v \in V$ we have
\begin{equation*}
R(u) \cdot R(v) = E_r^{-1}(\beta(r,G)f(u) \cdot f(v)),
\end{equation*}
and moreover the $R(u)$ are unit vectors. Hence, one can use the Cholesky decomposition of~$(R(u) \cdot R(v)) \in \R^{V \times V}_{\succeq 0}$ to define the desired function $g\colon V \to S^{|V|-1}$.
\end{proof}

We conclude this section with a few remarks on the lemma and its proof:

\begin{enumerate}
\item To approximate the Gram matrix $(R(u) \cdot R(v))$ it is enough to compute the series expansion of $E^{-1}_r$ and the matrix $U$ to the desired precision. The latter is found by solving a semidefinite program.
\item
Krivine proved the statement of the lemma in the case $r = 1$ and for bipartite graphs~$G$; then, $\vartheta(\overline{G}) = 2$ holds. Here, one only needs the first step of the proof. Also, $\beta(1,G)$ can be computed analytically. We have $E_1^{-1}(t) = \sin(\pi/2 t)$ and 
\begin{equation*}
\sum_{k=0}^{\infty} \left|(-1)^{2k+1} \frac{(\pi/2)^{2k+1}}{(2k+1)!}\right| t^{2k+1} = \sinh(\pi/2 t).
\end{equation*}
Hence, $\beta(1,G) = 2 \arcsinh(1) / \pi= 2 \ln(1+\sqrt{2}) / \pi$.
\item
In the second step one can also work with any feasible solution of the semidefinite program~\eqref{opt:theta-gbar}. For instance one can replace $\vartheta(\overline{G})$ in the lemma by the chromatic number~$\chi(G)$ albeit getting a potentially weaker bound.
\item
Alon, Makarychev, Makarychev, and Naor \cite{AlonMakarychevMakarychevNaor} also gave an upper bound for $K(1,G)$ using the theta number of the complement of~$G$. They prove that
\begin{equation*}
K(1,G) \leq O(\log\vartheta(\overline{G})),
\end{equation*}
which is much better than our result in the case of large $\vartheta(\overline{G})$. However, our bound is favourable when $\vartheta(\overline{G})$ is small.
In Section~\ref{sec:highchrom} we generalize the methods of Alon, Makarychev, Makarychev, and Naor~\cite{AlonMakarychevMakarychevNaor} to obtain better upper bounds on~$K(r,G)$ for $r\geq 2$ and large $\vartheta(\overline{G})$.
\item
Finally, notice that in the first step it was essential that the Taylor expansion of $E_r^{-1}$ has convergence radius of at least one.
\end{enumerate}

\section{A refined, dimension-dependent analysis}
\label{refined section}

So far we only compared the two problems $\sdp_{\infty}$ and $\sdp_r$. One can perform a refined, dimension-dependent analysis by comparing $\sdp_q$ and $\sdp_r$ when $q \geq r$. 

Let~$K(q \mapsto r, G)$, where~$q \geq r$, be the least number such that
\begin{equation*}
\sdp_q(G,A) \leq K(q \mapsto r, G) \sdp_r(G,A)
\end{equation*}
for all~$A\colon V \times V \to \R$. In this section we give an upper bound for~$K(q \mapsto r, G)$ that depends on~$q$ and~$r$. For fixed~$r$, this upper bound will become smaller as~$q$ comes closer to~$r$.  Krivine \cite{Krivine} gave such a refined, dimension-dependent analysis in the bipartite case; he showed that
\[
K(2 \mapsto 1, K_{n,n}) = \sqrt{2},\quad K(3 \mapsto 1, K_{n,n}) \leq 1.517, \quad\text{and}\quad K(4 \mapsto 1, K_{n,n}) \leq \pi/2.
\]
Avidor and Zwick~\cite{AvidorZwick} analyzed the cases $r=1$ and $q\in\{2,3\}$ for bipartite~$G=K_{n,n}$ and $A$ the Laplacian of a graph~$G'$ on~$n$ nodes.

Our upper bound comes from the following lemma:

\begin{lemma}
\label{lem:sphere-embedding}
Let~$G = (V, E)$ be a graph with at least one edge. Given~$f\colon V
\to S^{q-1}$, there is a function~$g\colon V \to S^{|V|-1}$ such
that whenever $u$ and $v$ are adjacent, then
\begin{equation*}
E_r(g(u) \cdot g(v)) = \beta(q \mapsto r,G) f(u) \cdot f(v),
\end{equation*}
where~$0 < \beta(q \mapsto r, G) \leq 1$  is such that~$\beta(q \mapsto r, G) > \beta(q + 1 \mapsto r, G)$ and~$\beta(q \mapsto r, G) > \beta(r,G)$ for all~$q \geq 2$.
\end{lemma}

\noindent The proof of the lemma will also give a procedure to compute $\beta(q \mapsto r, G)$ explicitly.

So we have the theorem:

\begin{theorem}
Let~$G = (V, E)$ be a graph with at least one edge and let~$q \geq r \geq 1$ be integers. Then~$K(q \mapsto r, G) \leq \beta(q \mapsto r, G)^{-1}$.
\end{theorem}

\begin{proof}
Combine Lemma~\ref{lem:sphere-embedding} with Algorithm~A from Section~\ref{sec:our methods}.
\end{proof}

The proof of the lemma uses some basic facts from harmonic
analysis, which we now summarize. For measurable functions~$f$,
$g\colon [-1, 1] \to \R$ we consider the inner product
\begin{equation}
\label{eq:inner-p}
\langle f, g \rangle_n = \int_{-1}^1 f(t) g(t) (1 - t^2)^{(n-3)/2}\,
dt.
\end{equation}
We say that a continuous function~$f\colon [-1, 1] \to \R$ is of
\textit{positive type for~$S^{n-1}$} if for any choice~$x_1$,
\dots,~$x_N$ of points in~$S^{n-1}$ we have that the
matrix~$\bigl(f(x_i \cdot x_j)\bigr)_{i, j = 1}^N$ is positive
semidefinite. If two continuous functions~$f$, $g\colon [-1, 1] \to
\R$ are of positive type for~$S^{n-1}$, then~$\langle f, g \rangle_n
\geq 0$.

Schoenberg~\cite{Schoenberg} characterized the continuous functions
of positive type in terms of Gegenbauer polynomials. We denote
by~$P_k^n$ the Gegenbauer polynomial of degree~$k$ and
parameter~$(n-2)/2$ which is normalized so that~$P_k^n(1) = 1$. Notice
that this normalization is not the one commonly found in the
literature.

The Gegenbauer polynomials~$P_0^n$, $P_1^n$, $P_2^n$, \dots\ are
pairwise orthogonal with respect to the inner
product~\eqref{eq:inner-p}, and they form a complete orthogonal system
for the space~$L^2([-1, 1])$, equipped with the inner
product~\eqref{eq:inner-p}. 

Schoenberg's characterization of the functions of positive type is as
follows: A function~$f\colon [-1, 1] \to \R$ is continuous and of positive type
for~$S^{n-1}$ if and only if
\begin{equation}
\label{eq:pos-decomp}
f(t) = \sum_{k=0}^\infty a_k P_k^n(t)
\end{equation}
for some nonnegative numbers~$a_0$, $a_1$, $a_2$, \dots\ such
that~$\sum_{k=0}^\infty a_k$ converges, in which case the series
in~\eqref{eq:pos-decomp} converges absolutely and uniformly in~$[-1,
1]$.

A continuous function~$f\colon [-1, 1] \to \R$ can also be of positive
type for spheres of every dimension. Schoenberg~\cite{Schoenberg}
also characterized these functions. They are the ones that can be
decomposed as
\begin{equation*}
f(t) = \sum_{k=0}^\infty a_k t^k
\end{equation*}
for some nonnegative numbers~$a_0$, $a_1$, $a_2$, \dots\ such
that~$\sum_{k=0}^\infty a_k$ converges.

A polynomial in~$\R[x_1, \ldots, x_n]$ is \textit{harmonic} if it is
homogeneous and vanishes under the Laplace operator~$\Delta
= \partial^2 / \partial x_1^2 + \cdots + \partial^2 / \partial
x_n^2$. Harmonic polynomials restricted to the unit sphere $\sphere{n}$ are related to Gegenbauer polynomials by
the \textit{addition formula} (see e.g.\ Andrews, Askey, and
Roy~\cite[Theorem 9.6.3]{AndrewsAskeyRoy}): Let~$H_k$ be the space of degree~$k$ harmonic polynomials on~$n$
variables. Any orthonormal basis of~$H_k$ can be scaled so as to give
a basis~$e_{k, 1}$, \dots,~$e_{k, h_k}$ of~$H_k$ for which the
following holds: For every~$u$, $v \in S^{n-1}$ we have that
\begin{equation*}
P_k^n(u \cdot v) = \sum_{i=1}^{h_k} e_{k,i}(u) e_{k,i}(v).
\end{equation*}

With this we have all that we need to prove the lemma. We only consider the bipartite case in the proof in order to simplify the notation and to make the argument more transparent. One can handle the nonbipartite case exactly in the same way as in the proof of Lemma~\ref{lem:gen-embedding}.

\begin{proof}[Proof of Lemma~\ref{lem:sphere-embedding}]
As before, we construct the function $g:V\to\sphere{|V|}$ from functions $S$ and $T$ that satisfy $S(x)\cdot T(y) = E^{-1}_r(\beta x\cdot y)$ for some real number $\beta$.

Fix~$0 < \beta \leq 1$ and consider the expansion
\begin{equation*}
E_r^{-1}(\beta t) = \sum_{k=0}^\infty g_k^q(\beta) P_k^q(t),
\end{equation*}
which converges in the~$L^2$ sense. 

\begin{claim}
The series~$\sum_{k=0}^\infty
|g_k^q(\beta)|$ converges for every~$0 < \beta \leq 1$, and hence the
above expansion converges absolutely and uniformly for~$t \in [-1,
1]$.
\end{claim}

\begin{claimproof}
To prove the claim we use the comparison test. To this end, consider the expansion~$E_r^{-1}(t) = \sum_{k=0}^\infty b_k t^k$ and recall that $\sum_{k=0}^\infty |b_k|$ converges. We may define the function~$\bE_r^{-1}(t) = \sum_{k=0}^\infty |b_k| t^k$, which is then of positive type for every sphere. So by Schoenberg's theorem we can write
\begin{equation*}
\bE_r^{-1}(t) = \sum_{k=0}^\infty \bg_k^q P_k^q(t)
\end{equation*}
for nonnegative numbers~$\bg_k^q$ such that~$\sum_{k=0}^\infty \bg_k^q$ converges. Now notice that
\begin{equation*}
g_k^q(\beta) = \|P_k^q\|^{-2}_q \langle E_r^{-1}(\beta t), P_k^q
\rangle_q = \|P_k^q\|^{-2}_q \sum_{l=0}^\infty b_l \beta^l \langle
t^l, P_k^q\rangle_q,
\end{equation*}
where~$\|P_k^q\|_q = \langle P_k^q, P_k^q \rangle_q^{1/2}$. Above,
since~$t^l$ is a function of positive type for every sphere, we 
have that~$\langle t^l, P_k^q \rangle_q \geq 0$. But we also have that
\begin{equation*}
\bg_k^q = \|P_k^q\|^{-2}_q \langle \bE_r^{-1}, P_k^q
\rangle_q = \|P_k^q\|^{-2}_q \sum_{l=0}^\infty |b_l| \langle
t^l, P_k^q\rangle_q,
\end{equation*}
and we see that~$|g_k^q(\beta)| \leq \bg_k^q$ for all~$k \geq 0$
and~$0 < \beta \leq 1$. This proves the claim.  
\end{claimproof}

From the proof of the claim, it it also clear that
\begin{equation}
\label{eq:sum-beta}
\sum_{k=0}^\infty |g_k^q(\beta)|
\end{equation}
is a continuous function of~$\beta$.

Now, let~$\beta(q \mapsto r, G)$ be the maximum number~$\beta\in(0, 1]$ that is such that
\begin{equation*}
\sum_{k=0}^\infty |g_k^q(\beta)| = 1.
\end{equation*}
By the intermediate value theorem, such a number exists because~\eqref{eq:sum-beta} is continuous as a
function of~$\beta$, being equal to~$0$ when~$\beta = 0$ and at
least~$E_r^{-1}(1) = 1$ when~$\beta = 1$.

Consider the Hilbert space
\begin{equation*}
\Hcal = \bigoplus_{k=0}^\infty \R^{h_k},
\end{equation*}
equipped with the Euclidean inner product and where $h_k$ is the dimension of $H_k$, the space of harmonic polynomials of degree~$k$ on~$q$ variables. For a vector~$x \in S^{q-1}$,
consider the vectors~$S(x)$ and~$T(x) \in \Hcal$ given componentwise by
\begin{equation*}
\begin{split}
S(x)_k&= \sqrt{|g_k^q(\beta(q \mapsto r, G))|} (e_{k, 1}(x), \ldots, e_{k,
  h_k}(x))\qquad\text{and}\\
T(x)_k&= \sign(g_k^q(\beta(q \mapsto r, G))) \sqrt{|g_k^q(\beta(q \mapsto r, G))|} (e_{k, 1}(x), \ldots, e_{k,
  h_k}(x)).
\end{split}
\end{equation*}
By the addition formula  we have that
\begin{equation*}
S(f(u)) \cdot T(f(v)) = E_r^{-1}(\beta(q \mapsto  r, G) f(u) \cdot f(v)).
\end{equation*}
Moreover, we also have that
\begin{equation*}
\|S(f(u))\|^2 = \|T(f(v))\|^2 = \sum_{k=0}^\infty |g_k^q(\beta(q \mapsto r, G))| = 1,
\end{equation*}
and so from the Gram matrix of the vectors~$S(f(u))$ and~$T(f(v))$ we may obtain the function $g\colon V \to S^{|V|-1}$ as we wanted.

Now we show that~$\beta(q \mapsto r, G) > \beta(q + 1 \mapsto r, G)$ for all~$q \geq 2$. To this
end, consider the function
\begin{equation*}
F_\beta(t) = \sum_{k=0}^\infty |g_k^{q+1}(\beta)| P_k^{q+1}(t).
\end{equation*}
Since the series~$\sum_{k=0}^\infty |g_k^{q+1}(\beta)|$ converges,
from Schoenberg's theorem we see that~$F_\beta$ is a continuous function
of positive type for the sphere~$S^q$. Notice moreover that, by
definition,~$\beta(q+1 \mapsto r, G)$ is the maximum number~$\beta\in(0, 1]$ such
that~$F_\beta(1) = 1$.

Since~$F_\beta$ is of positive type for~$S^q$, it is also of positive type
for~$S^{q-1}$, and then we may write
\begin{equation*}
F_\beta(t) = \sum_{k=0}^\infty a_k(\beta) P_k^q(t),
\end{equation*}
and we have~$\sum_{k=0}^\infty a_k(\beta) = 1$ as for all $k$, $P_k^q(1) = 1$. We also have
the expression
\begin{equation}
\label{eq:ak-exp}
a_k(\beta) = \|P_k^q\|^{-2}_q \langle F_\beta, P_k^q\rangle_q
           = \|P_k^q\|^{-2}_q \sum_{l=0}^\infty |g_l^{q+1}(\beta)|
           \langle P_l^{q+1}, P_k^q\rangle_q.
\end{equation}
Notice that, since both~$P_l^{q+1}$ and~$P_k^q$ are of positive type
for~$S^{q-1}$, we have that $\langle P_l^{q+1}, P_k^q\rangle_q \geq 0$
for all~$l$ and~$k$.

Now, from the expansion
\[
E_r^{-1}(\beta t) = \sum_{k=0}^\infty g_k^{q+1}(\beta) P_k^{q+1}(t)
\]
we see that
\begin{equation}
\label{eq:small-exp}
g_k^q(\beta) = \|P_k^q\|^{-2}_q \langle E_r^{-1}(\beta t),
P_k^q\rangle_q = \|P_k^q\|^{-2}_q \sum_{l=0}^\infty g_l^{q+1}(\beta)
\langle P_l^{q+1}, P_k^q\rangle_q.
\end{equation}
The function $E^{-1}_r$ is not of positive type because the coefficient $b_3$ of its Taylor expansion is always negative (this can easily be checked using Eq.~\eqref{eq:rankgroth-mfinverse}), and so
some of the~$g_l^{q+1}(\beta)$ must be negative. This, together
with~\eqref{eq:ak-exp} and~\eqref{eq:small-exp}, implies
that~$|g_k^q(\beta)| < a_k(\beta)$ for all~$0 < \beta \leq 1$. So we
must have that
\begin{equation*}
\sum_{k=0}^\infty |g_k^q(\beta(q+1 \mapsto r, G))| < \sum_{k=0}^\infty a_k(\beta(q+1 \mapsto
r, G)) = 1,
\end{equation*}
and we see that~$\beta(q \mapsto r, G) > \beta(q+1 \mapsto r, G)$.  In a similar way, one may show that~$\beta(q \mapsto r, G) >
\beta(r, G)$.
\end{proof}

\section{Better bounds for large chromatic numbers}
\label{sec:highchrom}

For graphs with large chromatic number, or more precisely with large $\vartheta(\overline G)$, our bounds on $K(r,G)$ proved above can be improved using the techniques of Alon, Makarychev, Makarychev, and Naor~\cite{AlonMakarychevMakarychevNaor}. In this section, we show how their bounds on $K(1,G)$ generalize to higher values of~$r$.

\begin{theorem}\label{thm:rankgroth-largechrom}
Given a graph $G = (V,E)$ and integer $1\leq r\leq \log\vartheta(\overline G)$, we have
\beqn
K(r,G) \leq O\left(\frac{\log\vartheta(\overline G)}{r}\right).\eeqn
\end{theorem}

\begin{proof}
It suffices to show that for any matrix $A:V\times V\to\R$, we have
\beqn
\sdp_r(G,A)\geq  \Omega\left(\frac{r}{\log\vartheta(\overline G)}\right)\sdp_\infty(G,A).
\eeqn

Fix a matrix $A:V\times V\to\R$. Let $f:V\to\sphere{|V|}$ be optimal for~$\sdp_{\infty}(G,A)$, so that
\beqn
\sum_{\{u,v\}\in E} A(u,v) f(u)\cdot f(v) = \sdp_{\infty}(G,A).
\eeqn

Let $\lambda = \vartheta(\overline G)$, and $\widetilde Z:V\times V\to\R$ be an optimal solution of~\eqref{opt:theta-gbar}. Let~$J$ be the $2|V| \times 2|V|$ all-ones matrix and $I$ the $2 \times 2$ identity matrix. Since the matrix $(I\otimes \widetilde Z + J)/\lambda$ is positive semidefinite, we  obtain from its Gram decomposition functions $s,t:V\to\R^{2|V|}$ that satisfy
\begin{enumerate}
\item $s(u)\cdot s(u) = t(u)\cdot t(u) = 1$ for all $u\in V$.
\item $s(u)\cdot t(u) = 1/\lambda$ for all $u,v\in V$.
\item $s(u)\cdot s(v) = t(u)\cdot t(v) = 0$ for all $\{u,v\}\in E$.
\end{enumerate}
\medskip

Let $\HS$ be the Hilbert space of vector-valued functions $h:\R^{r\times |V|}\to \R^r$ with inner product
\beqn
(g,h) = \Exp[g(Z)\cdot h(Z)],
\eeqn
where the expectation is taken over random $r \times |V|$ matrices $Z$ whose entries are i.i.d.~$N(0,1/r)$ random variables.

Let $R\geq 2$ be some real number to be set later.
Define for every $u\in V$ the function $g_u\in \HS$ by
\beqn
 g_u(Z) =\left\{\begin{array}{ll}
		\frac{Zf(u)}{R} & \text{if } \|Zf(u)\| \leq R\\[.2cm]
		\frac{Zf(u)}{\|Zf(u)\|} & \text{otherwise.}
		\end{array}
		\right.
\eeqn
Notice that for every matrix $Z\in \R^{r\times |V|}$, the vector $g_u(Z)\in\R^r$ has Euclidean norm at most 1. It follows by linearity of expectation that
\beqn
\sdp_r(G,A) \geq \Exp\bigg[\sum_{\{u,v\}\in E} A(u,v)\, g_u(Z)\cdot g_v(Z)\bigg] = \sum_{\{u,v\}\in E} A(u,v) (g_u,g_v).
\eeqn
We proceed by lower bounding the right-hand side of the above inequality.

Based on the definition of $g_u$ we define two functions $h_u^0,h_u^1\in \HS$ by
\begin{align*}
h_u^0(Z) = \frac{Z f(u)}{R} +  g_u(Z)& &\text{and}& &h_u^1(Z) &=& \frac{Z f(u)}{R} -  g_u(Z).
\end{align*}
For every $u\in V$, define the function $H_u\in \R^{2|V|}\otimes \HS$ by		
\beqn
H_u = \frac{1}{4} s(u)\otimes h_u^0 + 2\lambda\, t(u)\otimes h_u^1.
\eeqn
We  expand the inner products $(g_u,g_v)$ in terms of $f(u)\cdot f(v)$ and $\langle H_u,H_v\rangle$.

\begin{claim}
For every $\{u,v\}\in E$ we have
\beqrn
(g_u,g_v) &=& \frac{1}{R^2} f(u)\cdot f(v) - \langle H_u,H_v\rangle.
  \eeqrn
\end{claim}

\begin{claimproof}
Simply expanding the inner product $\langle H_u,H_v\rangle$ gives
\[
\begin{split}
\langle H_u,H_v\rangle &= \frac{s(u)\cdot s(v)}{16} (h_u^0,h_v^0)\, +\, 4\lambda^2\big(t(u)\cdot t(v)\big)\,(h_u^1,h_v^1) \\[.2cm]
&\qquad{} + \frac{\lambda}{2}\Big[ \big(s(u)\cdot t(v)\big)\, (h_u^0,h_v^1)\, +\, \big(t(u)\cdot s(v)\big)\, (h_u^1,h_v^0)\Big].
\end{split}
\]
It follows from property 3 of $s$ and $t$ that the above terms involving $s(u)\cdot s(v)$ and $t(u)\cdot t(v)$ vanish. By property 2, the remaining terms reduce to
\[
\begin{split}
\frac{1}{2}\Big((h_u^0,h_v^1) + (h_u^1,h_u^0)\Big) &=  \frac{1}{2}\Exp\left[\left(\frac{Zf(u)}{R} +  g_u(Z)\right)\cdot \left(\frac{Zf(v)}{R} -  g_v(Z)\right)\right]\\
 &\qquad{} + \frac{1}{2}\Exp\left[\left(\frac{Zf(u)}{R} -  g_u(Z)\right)\cdot \left(\frac{Zf(v)}{R} +  g_v(Z)\right)\right].
\end{split}
\]

 Expanding the first expectation gives
\beqn
 \frac{1}{R^2}\Exp[f(u)^{\sf T}Z^{\sf T}Zf(v)] - (g_u,g_v)-
 \Exp\left[\frac{Zf(u)}{R}\cdot  g_v(Z)\right] + 
  \Exp\left[ g_u(Z)\cdot \frac{Zf(v)}{R}\right]
\eeqn
and expanding the second gives
\beqn
 \frac{1}{R^2}\Exp[f(u)^{\sf T}Z^{\sf T}Zf(v)] - (g_u,g_v)+
 \Exp\left[\frac{Zf(u)}{R}\cdot  g_v(Z)\right] -
  \Exp\left[ g_u(Z)\cdot \frac{Zf(v)}{R}\right].
\eeqn
Adding these two gives that the last two terms cancel. Since $\Exp[Z^{\sf T}Z] = I$, what remains equals
\beqn
\frac{1}{R^2} f(u)\cdot f(v) - (g_u,g_v),
\eeqn
which proves the claim.
\end{claimproof}

From the above claim it follows that
\[
\begin{split}
\sum_{\{u,v\}\in E}A(u,v) (g_u,g_v)
&=  \frac{1}{R^2}\sdp_{\infty}(G,A) - \sum_{\{u,v\}\in E}A(u,v) \langle H_u,H_v\rangle\\
&\geq \left(\frac{1}{R^2} - \max_{u\in V}\| H_u\|^2\right) \sdp_\infty(G,A),
\end{split}
\]
where  $\|H_u\|^2 = \langle H_u,H_u\rangle$.

By the triangle inequality, we have for every $u\in V$,
\beqn
\|H_u\|^2 \leq \left(\frac{1}{4} \|h^0_u\| + 2 \lambda \|h^1_u\|\right)^2
\leq 
\frac{1}{R^2}\left(\frac{1}{2} + 2\lambda R\, \Exp\left[\Big\|\frac{ Z f(u)}{R} -  g_u(Z)\Big\|\right]\right)^2.
\eeqn

By the definition of $ g_u$, the vectors $Zf(u)$ and $ g_u$ are parallel. Moreover, they are equal if $\|Zf(u)\|\leq R$. Since $f(u)$ is a unit vector, the $r$ entries of the random vector $Zf(u)$ are i.i.d.~$N(0,1/r)$ random variables. Hence, 
\[
\begin{split}
\Exp\left[\Big\|\frac{ Z f(u)}{R} -  g_u(Z)\Big\|\right] &= \int_{\R^r}\ind[\|x\|\geq R]\Big(\frac{\|x\|}{R} - 1\Big)\Big(\frac{r}{2\pi}\Big)^{r/2} e^{-r\|x\|^2/2}dx\\
&= \int_R^{\infty} \int_{\sphere{r}}\rho^{r-1}\Big(\frac{\rho}{R} - 1\Big)\Big(\frac{r}{2\pi}\Big)^{r/2} e^{-r\rho^2/2}d\rho d\tilde\omega_r(\xi)\\
&\leq \frac{r^{r/2}}{R\Gamma\big(\frac{r}{2}\big)}\int_R^{\infty} \rho^r e^{-r\rho^2/2}d\rho,
\end{split}
\]
where $\tilde\omega_r$ is the unique rotationally invariant measure on $\sphere{r}$, normalized such that $\tilde\omega_r(\sphere{r}) = r^{r/2}/\Gamma(r/2)$.
Using a substitution of variables, we get
\beqn
\int_R^{\infty} \rho^r e^{-r\rho^2/2}d\rho = \frac{1}{2}\Big(\frac{2}{r}\Big)^{(r+1)/2} \Gamma\Big(\frac{r+1}{2},\frac{rR^2}{2}\Big),
\eeqn
where $\Gamma(a,x)$ is the lower incomplete Gamma function~\cite[Eq.~(4.4.5)]{AndrewsAskeyRoy}. 

Collecting the terms from above then gives the bound
\beq\label{eq:chrom5}
\sdp_r(G,A)\geq \frac{1}{R^2}\left(1 -\left(\frac{1}{2} +  \lambda\frac{2^{(r+1)/2}}{\sqrt{r}\Gamma\big(\frac{r}{2}\big)}\Gamma\Big(\frac{r+1}{2}, \frac{rR^2}{2}\Big)\right)^2\right)\sdp_{\infty}(G,A).
\eeq
The bound in the theorem follows by setting $R$ as small as possible such that the above factor between brackets is some positive constant.

By Stirling's formula, there is a constant $C_1>0$ such that $\Gamma(x) \geq C_1 e^{-x}x^{x-1/2}$ holds (see for example Abramowitz and Stegun~\cite[Eq.~(6.1.37)]{Abramowitz:1964}). Hence, for some constants $c,C>0$, we have
\beq\label{eq:chrom2}
\frac{2^{(r+1)/2}}{\sqrt{r}\Gamma\big(\frac{r}{2}\big)} \leq  C\left(\frac{c}{r}\right)^{r/2}
\eeq

The power series  of the incomplete gamma function~\cite[Eq.~(6.5.32)]{Abramowitz:1964} gives that if $a\leq x$, for some constant $C_2>0$, the inequality $\Gamma(a,x) \leq C_2x^ae^{-x}$ holds.
As $R \geq 2$, for some constants $d,D>0$, we have
\beq\label{eq:chrom3}
\Gamma\left(\frac{r+1}{2}, \frac{rR^2}{2}\right) \leq D \sqrt{r} \left(\frac{r}{d^{R^2}}\right)^{r/2}.
\eeq
Putting together estimates~\eqref{eq:chrom2} and~\eqref{eq:chrom3} gives
\beqrn
\lambda\frac{2^{(r+1)/2}}{\sqrt{r}\Gamma\big(\frac{r}{2}\big)} \Gamma\left(\frac{r+1}{2}, \frac{rR^2}{2}\right) &\leq & CD \sqrt{r}\lambda\left(\frac{c}{d^{R^2}}\right)^{r/2}.
\eeqrn
Since $r\leq\log\lambda$ there is some constant $C'$ such that for $R^2 = C'\big(\log \lambda\big)/r$, the above value is less than $1/4$. It follows that for this value of $R$, Inequality~\eqref{eq:chrom5} is nontrivial and we get the result.
\end{proof}

%

%
%
%

\section*{Acknowledgements}

The third author thanks Assaf Naor for helpful comments, and the Institute for Pure \& Applied Mathematics at UCLA for its hospitality and support during the program ``Modern Trends in Optimization and Its Application'', September 17--December 17, 2010.

\end{document}